\documentstyle[12pt]{article}
\def\mineappendix{
        \setcounter{section}{1}
        \setcounter{subsection}{0}
        \def\thesection{\Alph{section}}
        \def\sectionap{\@startsection  {section}{1}{\z@}
                        {-3.5ex plus-1ex minus-.2ex} {0ex plus.2ex}
                        {\reset@font\Large\bf  Appendix:  \, }
                        }
        }
\makeatother
\def\Proclaim #1. #2\par{\bigbreak\noindent{\sc#1.\enspace}{\it#2}\par}

\font\Bbbfont=msbm10
\newfam\msbfam
\textfont\msbfam=\Bbbfont \textfont\msbfam=\Bbbfont

\def\pa{\partial}
\def\al{\alpha}

\def\si{\sigma}

\def\wt{\widetilde}

\def\la{\lambda}

\def\wt{\widetilde}

\def\nn{\nonumber}

\def\ga{\gamma}

\def\nd{\noindent}

\newtheorem{Theorem}{Theorem}
\newtheorem{Lemma}{Lemma}

\newtheorem{Proposition} {Proposition}
\newtheorem{Remark}{Remark}

\title{ Explicit Blowing-up Solutions to the Schr\"odinger Maps from ${\bf R}^2$ to the
Hyperbolic 2-Space ${\bf H}^2$}
\author{Qing Ding\\
Institute of Mathematics and Key Lab. of Math. for Nonlinear
Sciences\\ Fudan University, Shanghai
200433, China\\
 E-mail address: qding@fudan.edu.cn}
\date{}

\begin{document}
\maketitle

\begin{abstract}
In this article, we prove that the equation of the Schr\"odinger
maps from ${\bf R}^2$ to the hyperbolic 2-space ${\bf H}^2$ is
$SU(1,1)$-gauge equivalent to the following $1+2$ dimensional
nonlinear Schr\"odinger-type system of unknown three complex
functions $p,q,r$ and a real function $u$:
\begin{eqnarray}
\left\{\begin{array}{c} iq_t+q_{z{\bar z}}-2u
q+2({\bar p}q)_z-2pq_{\bar z}-4|p|^2q=0\qquad\\
ir_t-r_{z{\bar z}}+2u r+2({\bar p}r)_z-2pr_{\bar z}+4|p|^2r=0\qquad\\
ip_t+(qr)_{\bar z}-u_z=0\qquad\qquad\qquad\qquad\qquad\qquad~
\\{\bar p}_z+p_{\bar z}= -|q|^2+|r|^2,~-{\bar r}_z+q_{\bar
z}=-2(p{\bar r}+{\bar p}q),
\end{array} \right.   \nn
\end{eqnarray}
where $z$ is a complex coordinate of the plane ${\bf R}^2$ and $\bar
z$ is the complex conjugate of $z$. Though this nonlinear
Schr\"odinger-type system looks much complicated, it admits a class
of the following explicit blowing-up smooth solutions:
\begin{eqnarray}
p=0,~ q=={{e^{i{{bz\bar z}\over{2(a+bt)}}}}\over{a+bt}} \al \bar z,~
r= {{e^{-i{{bz\bar z}\over{2(a+bt)}}}}\over{(a+bt)}}\al \bar z,~
u={{2\al^2z{\bar z} }\over{(a+bt)^2}},\nn
\end{eqnarray}where $a,b$ are real numbers with $ab<0$ and $\al$ satisfies
$\al^2=\frac{b^2}{16}$. From these facts, we explicitly construct
smooth solutions to the Schr\"odinger maps from ${\bf R}^2$ to the
hyperbolic 2-space ${\bf H}^2$ by using the gauge transformations,
such that the absolute of their gradients blows up in finite time.
This reveals some blow-up phenomenon of Schr\"odinger maps.

\bigskip
\nd PACS numbers: 02.40.Ky; 02.60.Lj; 07.55.Db

\nd Keywords: Schr\"odinger maps, gauge transformation, blowing-up
solutions
\end{abstract}

\section *{\S 1. Introduction}
The Landau-Lifshitz equations, or in other words, the generalized
Heisenberg models for a continuous ferromagnetic spin vector
$s=(s_1,s_2,s_3)$ $\in S^2\hookrightarrow {\bf R}^3$ (see, for
example, \cite{LaL,PaT,SuSu}),
\begin{eqnarray}
{\bf s}_t={\bf s}{\times}\Delta_{{\bf R}^n}{\bf  s}, \quad x\in {\bf
R}^n,\quad n=1,2,3,\cdots, \label{0}
\end{eqnarray}
where $\times$ denotes the cross product in the Euclidean 3-space
${\bf R}^3$ and $\Delta_{{\bf R}^n}$ is the Laplacian operator on
${\bf R}^n$, are important equations appeared in magnetic fields
theory. Eqs.(\ref{0}) have been generalized to various Hermitian
symmetric spaces (see, for example \cite{TU}). The following is one
of them
\begin{eqnarray}
{\bf s}_t={\bf s}{\dot\times}\Delta_{{\bf R}^n} {\bf s}, \quad x\in
{\bf R}^n, ~~n=1,2,3,\cdots, \label{1}
\end{eqnarray}
where ${\bf s}=(s_1,s_2,s_3)\in {\bf H}^2=\{(s_1,s_2,s_3)\in{\bf
R}^{2+1} :|s|^2=s_1^2+s_2^2-s_3^2=-1, s_3<0\} \hookrightarrow {\bf
R}^{2+1}$ is a `unit' spin vector in the Minkowski 3-space ${\bf
R}^{2+1}$, and $\dot\times$ denotes the pseudo-cross product in
${\bf R}^{2+1}$, i.e. for ${\bf a}=(a_1,a_2,a_3), {\bf
b}=(b_1,b_2,b_3)\in {\bf R}^{2+1}$, ${\bf a}\dot\times {\bf
b}=(a_2b_3-a_3b_2,a_3b_1-a_1b_3,-(a_1b_2-a_2b_1))$. Eqs.(\ref{1})
are regarded as dual versions of Eqs.(\ref{0}) (see \cite{GS,NSU}).
In fact, Eqs.(\ref{0}) relate to the $SU(2)$ ferromagnetic spin
models and, meanwhile, Eqs.(\ref{1}) relate to the $SU(1,1)$
ferromagnetic spin models in physics (see, for example \cite{OP}).
The significance of studying Eq.(\ref{1}) with $n=2$ is also
mentioned in \cite{NSU}. We would like to point out that there are
two components of the surface $s_1^2+s_2^2-s_3^2=-1$ in ${\bf
R}^{2+1}$, i.e., $s_3>0$ and $s_3<0$. Without loss of generality, we
shall fix ${\bf H}^2$ as the above surface in ${\bf R}^{2+1}$ with
the component of $s_3>0$, unless otherwise stated.

The Landau-Lifshitz equations (\ref{0}) and their dual versions
(\ref{1}) are special cases and typical examples of the so-called
Schr\"odinger maps (\cite{CSU}, \cite{GS}, \cite{NSU}) or
Schr\"odinger flows (\cite{DW2}, \cite{d1}) in geometry. A
Schr\"odinger map $u$ from a Riemannian manifold $(M,g)$ to a
K\"ahler manifold $(N,J)$ is defined to be a solution to the
(infinite dimensional) Hamiltonian system of the energy function
$E(u)=\int_M|\nabla u|^2dv_g$ on the mapping space $C^k(M,N)$ for
some $k>0$. More explicitly, the equation of the  Schr\"odinger
maps $u:M \to N$ is expressed as the following evolution system:
$$\begin{array}{l}
u_t=J(u)\tau(u),
\end{array}\label{2}$$
where $J$ is the complex structure on $N$ and $\tau(u)$ is the
tension field of $u$. It is a straightforward verification that the
Schr\"odinger  map from the Euclidean $n$-space ${\bf R}^n$ to the
2-sphere $S^2\hookrightarrow{\bf R}^3$ is just Eq.(\ref{0}) (for
example see \cite{CSU}, \cite{DW2}) and the Schr\"odinger map from
the Euclidean $n$-space ${\bf R}^n$ to the hyperbolic 2-space ${\bf
H}^2\hookrightarrow{\bf R}^{2+1}$ is nothing else than Eq.(\ref{1})
(\cite{d1}).

Though Eqs.(\ref{0}) (resp. Eqs.(\ref{1}))  have a unified
expressions for $n\ge1$, there are great differences between
dynamical properties of Eq.(\ref{0}) (resp. Eq.(\ref{1}))  with
$n=1$ and those of Eqs.(\ref{0}) (resp. Eqs.(\ref{1}))  with
$n\ge2$. When $n=1$, Eq.(\ref{0}) and (\ref{1}) are integrable and
can be solved by the method of inverse scattering techniques
(\cite{FaT}). Moreover, it was proved by Zarkharov and Takhtajan in
\cite{ZaT} that Eq.(\ref{0}) is gauge equivalent to the integrable
focusing nonlinear Schr\"odinger equation: $i\phi_t+\phi_{xx}+2
|\phi|^2\phi=0$ and, meanwhile, it was showed in \cite{d1} that
Eq.(\ref{1}) with $n=1$ is gauge equivalent to the defocusing
nonlinear Schr\"odinger equation: $q_t+q_{xx}-2|q|^2q=0$. When
$n\ge2$, Eqs.(\ref{0}) or Eqs.(\ref{1}) are, roughly speaking,
non-integrable and the understanding of their dynamical properties
becomes much more difficult than that of their one dimensional
cases. The local in time existence and uniqueness of
$W^{m,\sigma}({\bf R}^n)$-solutions to the Cauchy problem of the
Schr\"odinger map from ${\bf R}^n$ ($n\ge2$) to 2-sphere
$S^2\hookrightarrow R^3$ was established by Sulem, Sulem and Bardos
in 1986 in \cite{SuSu}, where $\sigma\ge 2$. They also proved the
global in time $W^{m+1,6}({\bf R}^n)$-existence of the Cauchy
problem of Eq.(\ref{0}) for small initial data. In 2000, Chang,
Shatah and Uhlenbeck displayed in \cite{CSU} the global in time
$W^{2,4}({\bf R}^2)$-existence of radially symmetric solutions to
the Cauchy problem of Schr\"odinger maps from ${\bf R}^2$ to
Riemannian surfaces for small initial data.  There is also previous
work \cite{DW2} by W.Y. Ding and Wang in local time existence and
uniqueness of solutions for Schr\"odinger maps to K\"ahler
manifolds. Grillakis and Stefanopoulos displayed in \cite{GS}
conservation laws and localized energy estimates of Schr\"odinger
maps to Riemannian surfaces. Recently, Nahmod, Stefanov and
Uhlenbeck proved the local well-posedness of the Cauchy problem of
the Schr\"odinger maps from ${\bf R}^2$ to the 2-sphere ${\bf S}^2$
and the hyperbolic 2-space ${\bf H}^2$ by using their equivalent
equations which are so-called the modified Schr\"odinger map
equations in \cite{NSU}. However, it is widely believed that there
should exist blowing-up solutions to Schr\"odinger maps (or flows)
(see, for example, \cite{Di}) when the dimension of the starting
manifolds is greater than 1.

In this paper, by using the geometric concept of gauge equivalence
for PDEs with prescribed curvature representation applied in \cite
{dz,dl}, we show that the equation of the Schr\"odinger maps from
${\bf R}^2$ to the hyperbolic 2-space ${\bf H}^2 \hookrightarrow
{\bf R}^{2+1}$: ${\bf s}_t={\bf
s}\dot{\times}\left({{\partial^2}\over{\partial x_1^2}}+
{{\partial^2}\over{\partial x_2^2}}\right){\bf s}$ is gauge
equivalent to the following $1+2$ dimensional nonlinear
Schr\"odinger-type system:
\begin{eqnarray}
\left\{\begin{array}{c} iq_t+q_{z{\bar z}}-2u
q+2({\bar p}q)_z-2pq_{\bar z}-4|p|^2q=0\\
ir_t-r_{z{\bar z}}+2u r+2({\bar p}r)_z-2pr_{\bar z}+4|p|^2r=0\\
ip_t+(qr)_{\bar z}-u_z=0\qquad\qquad\qquad\qquad\qquad
\end{array} \right.     \label{r8}
\end{eqnarray}
where $z=\frac{x_1+ix_2}{2}$ is a formal complex version of the
variables $x_1$ and $x_2$ ($\bar z$ is the complex conjugate of
$z$), $u$ is a unknown real function and $p,q,r$ are unknown complex
functions satisfying the following additional restrictions
\begin{eqnarray}
{\bar p}_z+p_{\bar z}= -|q|^2+|r|^2,\quad -{\bar r}_z+q_{\bar
z}=-2(p{\bar r}+{\bar p}q).  \label{rq=r}
\end{eqnarray}
This nonlinear Schr\"odinger-type system (\ref{r8},\ref{rq=r}) is
somewhat different from the modified Schr\"odinger map equations
deduced by Nahmod, Stefanov and Uhlenbeck in \cite{NSU}. By fixing
an ansatz solution: $p=0, q=e^{-i\theta}Q(\rho,t),
r=e^{-i\theta}{\bar Q}(\rho,t),
u=|Q|^2-2\int_{\rho}^{\infty}{{|Q|^2(\tau,t)}\over{\tau}}d\tau$ for
some functions $Q(t,\rho)$ to the system (\ref{r8},\ref{rq=r}),
where $(\rho,\theta)$ are the polar coordinates of $(x_1,x_2)$,
(\ref{r8}) leads $Q(t,\rho)$ to satisfy
\begin{eqnarray}
 iQ_t+\left(Q_{\rho\rho}+{1\over{\rho}}Q_{\rho}
-{1\over{\rho^2}}Q\right)-2Q\left(|Q|^2 -2\int_{\rho}^{\infty}
{{|Q(\tau,t)|^2}\over{\tau}}d\tau \right)=0.\label{QQ}
\end{eqnarray}
This is a well-known nonlinear Schr\"odinger-type equation gauged to
the radially symmetric Sch\"odinger map from ${\bf R}^2$ to the
hyperbolic 2-space ${\bf H}^2$ (see \cite{ CSU,La}). As
applications, by using the gauged nonlinear Schr\"odinger-type
system (\ref{r8},\ref{rq=r}), we first prove that there are no
nontrivial smooth radially symmetric harmonic maps from ${\bf R}^2$
to the Hyperbolic 2-Space ${\bf H}^2$. Then we show the existence of
class of smooth blowing-up solutions to the Schr\"odinger maps from
${\bf R}^2$ to the hyperbolic 2-space ${\bf H}^2$ by constructing
explicit blowing-up solutions to the nonlinear Schr\"odinger-type
system (\ref{r8},\ref{rq=r}). The latter reveals the blow-up
phenomenon of Schr\"odinger maps and gives an affirmative answer to
a problem proposed by W.Y. Ding in \cite{Di} for Schr\"odinger
flows. This also indicates that the introduction of the gauged
nonlinear Schr\"odinger-type system (\ref{r8},\ref{rq=r}) is much
more effective than that of the known Eq.(\ref{QQ}) in the study of
the Schr\"odinger maps from ${\bf R}^2$ to ${\bf H}^2$.

This paper is organized as follows. In the section 2, we  shall
transform the equation of the Schr\"odinger maps from ${\bf R}^2$ to
the hyperbolic 2-space ${\bf H}^2$ to its equivalent nonlinear
Schr\"odinger-type system (\ref{r8},\ref{rq=r}) by applying
geometric gauge theory. In section 3, as one of applications, we
show the existence of class of smooth blowing-up solutions to the
Schr\"odinger maps from ${\bf R}^2$ to the hyperbolic 2-space ${\bf
H}^2$ via the explicit construction of blowing-up solutions to the
gauged nonlinear Schr\"odinger-type system (\ref{r8},\ref{rq=r}).
Finally, in section 4, we close the paper with some conclusions and
remarks. We set an appendix at the end of the paper to give a
detailed proof of explicit general solutions to a matrix equation
which is important but not proved in the context of the paper.

\section * {\S 2. Gauge Equivalence}

Zakharov and Takhtajan introduced in \cite{ZaT} the geometric
concept of gauge equivalence between two soliton equations with zero
curvature representation, which provides a useful tool in the study
of integrable equations \cite{FaT}. In  \cite{dz,dl} the geometric
concept of gauge equivalence between (integrable) differential
equations with zero curvature representation has been generalized to
(nonintegrable or integrable) differential equations with prescribed
curvature representation. Now we find that this geometric idea is
applicable to the present situation, as we shall see below.

It is obvious that the Schr\"odinger maps from ${\bf R}^2$ to ${\bf
H}^2 \hookrightarrow {\bf R}^{2+1}$ (refer to (\ref{1}) with $n=2$)
is equivalent to the following system:
\begin{eqnarray}
{\bf s}_t=-{\bf s}{\dot\times}\left({{\partial^2}\over{\partial
x^2_1}}+ {{\partial^2}\over{\partial x_2^2}}\right){\bf s}
\label{r4}
\end{eqnarray}
by ${\bf s}\to -{\bf s}$. Notice that, in this case, the surface
${\bf H}=\{s_1^2+s_2^2-s_3^2=-1\}$ in ${\bf R}^{2+1}$ with $s_3>0$
becomes now the surface with $s_3<0$. Eq.(\ref{r4}) reads
explicitly:
\begin{eqnarray}
(s_1)_t&=&s_3(s_2)_{z\bar z}-s_2(s_3)_{z\bar z}\nn\\
(s_2)_t&=&s_1(s_3)_{z\bar z}-s_3(s_1)_{z\bar z}\nn \\
(s_3)_t&=&s_1(s_2)_{z\bar z}-s_2(s_1)_{z\bar z},\nn
\end{eqnarray}
where $z={{x_1+ix_2}\over2}$, ${\bar z}={{x_1-ix_2}\over2}$ which is
different from the usual complex version of the variables
$(x_1,x_2)$. Now we convert Eq.(\ref{r4}) to its matrix form:
\begin{eqnarray}
S_t=-{1\over2}[S, S_{z{\bar z}}], \label{r5}
\end{eqnarray}
where $S=\left(\begin{array}{cc}
       is_3 &s _1+is_2\\
       s_1-is_2 & -is_3
       \end{array}\right)\in su(1,1)$
        with $S^2=-I$.
Let's set
\begin{eqnarray}
A= Vd{\bar z} +i\lambda S d{z} +i\lambda\bigg(2V +S_{\bar z}S+2\al
S\bigg)dt \label{r6}
\end{eqnarray}
where $\lambda$ is a spectral parameter which is independent of
$t$,$z$ and $\bar z$, $V=V(\lambda,\xi,\eta,t)$ is a
$2\times2$-matrix function satisfying the equation:
\begin{eqnarray}
 (i\lambda S)_{\bar z}-V_{z} +[i\lambda S,V]=0 \label{rV}
\end{eqnarray}
and $\al=-{i\over{2}}\hbox{tr}(G^{-1}G_{\bar z}\sigma_3)$ is a
function depending only on $S$, here $G$ is an $SU(1,1)$-matrix
defined by (\ref{gamma}) below. $d+A$ can be geometrically
interpreted as defining a connection on the trivial principal bundle
${\bf R}^3\times SU(1,1)$ over ${\bf R}^3$ (the space of the
independent variables $x$, $y$ and $t$). From Yang-Mills theory, it
is a straightforward computation that the curvature $F_A$ of the
connection (\ref{r6}) is
\begin{eqnarray}
F_A&=&dA-A\wedge A\nn\\&=&\bigg\{-V_z+i\la S_{\bar z}-[V,i\la
S]\bigg\}d\bar z\wedge dz\nn\\&&+\bigg\{-V_t+i\lambda (2V +S_{\bar
z}S+2\al S)_{\bar z}-[V,i\lambda(2V
+S_{\bar z}S+2\al S)] \bigg\}d{\bar z} \wedge dt \nonumber \\
&&+\bigg\{-i\la S_t+i\la (2V +S_{\bar z}S+2\al S)_z-[i\la S,i\la (2V
+S_{\bar z}S+2\al S)]\bigg\} d{z} \wedge dt.\nn \label{cFA}
\end{eqnarray}
Thus the choice of $V$ in (\ref{rV}) above is to let the coefficient
of $d\bar z\wedge dz$ be vanish. Furthermore, set
\begin{eqnarray}
K&=&\bigg\{-V_t+i\lambda (2V +S_{\bar z}S+2\al S)_{\bar
z}-[V,i\lambda(2V
+S_{\bar z}S+2\al S)] \bigg\}d{\bar z} \wedge dt \nonumber \\
&&+i\lambda\bigg({1\over2}[S_{\bar z},S_{z}]+2(\al S)_{z}\bigg) d{z}
\wedge dt, \label{rK}
\end{eqnarray}
then it is also a direct calculation that the following prescribed
curvature condition:
\begin{eqnarray}
F_A=dA-A\wedge A=K \label{r7}
\end{eqnarray}
is actually equivalent to the matrix equation (\ref{r5}). Here in
the calculation we have used the fact $S^2=-I$ and hence the
identities: $S_{\bar z}S=-SS_{\bar z}$, $S_{\bar z z}S+S_{\bar
z}S_z+S_zS_{\bar z}=-SS_{\bar z z}$ by taking derivatives.

In this section, our aim is to use the prescribed curvature
representation formula (\ref{r7}) to transform the matrix equation
(\ref{r5}) into the (1+2) dimensional nonlinear Schr\"odinger-type
system (\ref{r8},\ref{rq=r}) in the category of Yang-Mills theory.

\begin{Theorem} For any given solution $S(t,z,{\bar z})$ to the
Schr\"odinger map from ${\bf R}^2$ to ${\bf H}^2$ (\ref{r5}), there
is a matrix-valued function $G(t,z,{\bar z})\in SU(1,1)$ such that
$S$ is transformed to a solution $(p(t,z,{\bar z}),q(t,z,{\bar
z}),r(t,z,{\bar z}),u(t,z,{\bar z}))$ to the 1+2 dimensional
nonlinear Schr\"odinger-type system (\ref{r8},\ref{rq=r}) by the
gauge transformation with $G$.
\end{Theorem}
\nd{\bf Proof}. Let $S=S(t,z,{\bar z})$ be a solution to
Eq.(\ref{r5}). We come to choose an $SU(1,1)$ matrix $G(t,z,{\bar
z})$ such that
\begin{eqnarray}
\sigma_3=iG^{-1}SG,\quad G^{-1}G_{z}=-\left(\begin{array}{cc}
  p & q\\
  {r}& -p
  \end{array}\right):=-U
  \label{14}
\end{eqnarray}
for some complex functions $p=p(t,z,{\bar z})$ and  $q=q(t,z,{\bar
z})$ and $r=r(t,z,{\bar z})$, where
$\sigma_3=\left(\begin{array}{cc}
       1& 0\\
       0& -1
       \end{array}\right)$ is the Pauli matrix. It is well-known
that the general $SU(1,1)$-solutions of the matrix equation
$\sigma_3=iG^{-1}SG$ are of the from (see the proof in the appendix
at the end of this paper):
\begin{eqnarray}
G={1\over{\sqrt{2(1-s_3)}}}(S-i\sigma_3){\rm
diag}(\gamma,{\bar\gamma}), \label{gamma}
\end{eqnarray}
where $\gamma$  is a complex function of $z$, $\bar z$ and $t$ with
$|\gamma|=1$.  For such a $SU(1,1)$-matrix $G$ given in
(\ref{gamma}), we have
\begin{eqnarray}
G_{x_1}=-G\left(\begin{array}{cc}
  is & \psi\\
  {\bar \psi}& -is
  \end{array}\right), \quad G_{x_2}=-G\left(\begin{array}{cc}
  il & \phi\\
  {\bar \phi}& -il
  \end{array}\right) \label{u}
\end{eqnarray}
for some real functions $s,l$ and complex functions $\phi, \psi$.
Since for the complex variables $z=(x_1+ix_2)/2$ and ${\bar
z}=(x_1-ix_2)/2$, we have $\frac{\pa}{\pa z}=\frac{\pa}{\pa
x_1}-i\frac{\pa}{\pa x_2}$ and $\frac{\pa}{\pa \bar
z}=\frac{\pa}{\pa x_1}+i\frac{\pa}{\pa x_2}$. Hence, from (\ref{u}),
\begin{eqnarray}
G_z=-G\left(\begin{array}{cc}
       l+is & (\psi-i \phi)\\
       ({\bar \psi}-i{\bar \phi})& -(l+is)
       \end{array}\right), \quad
G_{\bar z}=-G\left(\begin{array}{cc}
       -l+is & \psi+i \phi\\
       ({\bar \psi}+i{\bar \phi})& l-is
       \end{array}\right).       \nn
\end{eqnarray}
We thus have (\ref{14}) with $p=l+is$, $q=\psi-i\phi$ and $r={\bar
\psi}-i{\bar \phi}$, and
\begin{eqnarray}
G_{\bar z}=GP:=G \left(\begin{array}{cc}
       {\bar p} & -{\bar r}\\
       -{\bar q}& -{\bar p}
       \end{array}\right). \label{GP}
\end{eqnarray}
The integrability condition $P_{z}+U_{\bar z}+[P,U]=0$ of the linear
system: $G_{z}=-GU, G_{\bar z}=GP$ implies
\begin{eqnarray}
{\bar p}_z+p_{\bar z}= |r|^2-|q|^2,\quad -{\bar r}_z+q_{\bar
z}=-2(p{\bar r}+{\bar p}q), \nn
\end{eqnarray}
which are exactly the restrictions (\ref{rq=r}). Furthermore, we
have $\al=-i{\bar p}$ from (\ref{GP}) in the definition of the
connection $A$ given by (\ref{r6}). Now, for the connection $A$
given in (\ref{r6}) with $S$ being fixed above,  we  make the
following gauge transformation:
\begin{eqnarray}
A^{G}&=&-G^{-1}dG+G^{-1}A G\nn\\
&=&-G^{-1}(G_{\bar z}d{\bar z}+G_zdz+G_tdt)+G^{-1}\left(Vd{\bar z}
+i\lambda S d{z} +i\lambda(2V +S_{\bar z}S+2\al S)dt\right)G.\nn\\
   \label{AG1}
\end{eqnarray}
We would like to determine this $A^G$ explicitly and then prove the
theorem in this process. In fact, from the relations
(\ref{14},\ref{GP}) we have $S=-iG\sigma_3G^{-1}$, $S_{\bar
z}=-2iGP^{\hbox{(off-diag)}}$ $\si_3 G^{-1}$ and also $\al=-i{\bar
p}$. Hence
\begin{eqnarray}
S_{\bar z}S+2{\al} S
&=&-2iGP^{(\hbox{off-diag})}\si_3G^{-1}(-iG\si_3G^{-1})-2i\bar p(-iG\si_3G^{-1})\nn\\
&=&-2G[P^{\hbox{(off-diag)}}+\bar p\si_3]G^{-1} = -2G_{\bar
z}G^{-1}. \label{GG}
\end{eqnarray}
Thus substituting it into (\ref{AG1}) and using $G_z=-GU$, we obtain
\begin{eqnarray}
A^{G} &=& (-G^{-1}G_{\bar z}+G^{-1}VG)d{\bar
z}+(\lambda\sigma_3+U)d{z}\nn\\
&&+\bigg(-G^{-1}G_t+ 2i\lambda(G^{-1}VG-G^{-1}G_{\bar z})\bigg)dt.
   \label{AG}
\end{eqnarray}
As $A$ satisfying the prescribed curvature condition:
$F_A=dA-A\wedge A=K$, where $K$ is given by (\ref{rK}), from the
gauge theory we know that $A^G$ must fulfill
\begin{eqnarray}
F_{A^{G}}=dA^{G}-A^{G}\wedge A^{G}=G^{-1}(dA-A\wedge A)G=G^{-1}KG.
\label{017}
\end{eqnarray}
From (\ref{AG}) and by a direction computation we have
\begin{eqnarray}
F_{A^G}&=&\bigg\{(G^{-1}G_{\bar z}-G^{-1}VG)_z+(\la \si_3+U)_{\bar
z}-[-G^{-1}G_{\bar z}+G^{-1}VG,\la \si_3+U]\bigg\}d\bar z\wedge
dz\nn\\&&+\bigg\{-(-G^{-1}G_{\bar z}+G^{-1}VG)_t+(-G^{-1}G_t+2i\la
(G^{-1}VG-G^{-1}G_{\bar z}))_{\bar z}\nn\\&&-\bigg[(-G^{-1}G_{\bar
z}+G^{-1}VG),-G^{-1}G_t+2i\la
(G^{-1}VG-G^{-1}G_{\bar z})\bigg] \bigg\}d{\bar z} \wedge dt \nonumber \\
&&+\bigg\{-(\la \si_3+U)_t+(-G^{-1}G_t+2i\la (G^{-1}VG-G^{-1}G_{\bar
z}))_z\nn\\&&-\bigg[\la \si_3+U,-G^{-1}G_t+2i\la
(G^{-1}VG-G^{-1}G_{\bar z})\bigg]\bigg\} d{z} \wedge dt.\nn
\label{cFA}
\end{eqnarray}
Thus comparing respectively the coefficients of $d{z}\wedge d{\bar
z}$, $d{\bar z}\wedge dt$ and $d{z}\wedge dt$ in the both sides of
(\ref{017}), we have the following identities:
\begin{eqnarray}
-\bigg(-G^{-1}G_{\bar z}+G^{-1}VG\bigg)_{z}+U_{\bar
z}-\bigg[-G^{-1}G_{\bar z}+G^{-1}VG,\lambda \sigma_3+U\bigg]=0,
\label{FA1}
\end{eqnarray}
\begin{eqnarray}
&&-\bigg(-G^{-1}G_{\bar z}+G^{-1}VG\bigg)_t+\bigg(-G^{-1}G_t+
i\lambda(2G^{-1}VG-2G^{-1}G_{\bar z})\bigg)_{\bar z}
\nn \\
&&-\bigg[-G^{-1}G_{\bar z}+G^{-1}VG,-G^{-1}G_t+
i\lambda(2G^{-1}VG-2G^{-1}G_{\bar z})\bigg]
\nn\\
&=&G^{-1}\bigg\{-V_t+i\lambda (V_{\bar z}+S_{\bar z}S+2\al S)_{\bar
z}-\bigg[V,i\lambda (V_{\bar z}+S_{\bar z}S+2\al S)\bigg]\bigg\}G,
\label{FA2}
\end{eqnarray}
and
\begin{eqnarray}
&&-U_t+\bigg(-G^{-1}G_t+ 2i\lambda(G^{-1}VG-G^{-1}G_{\bar
z})\bigg)_{z}\nn
\\&& -\bigg[\lambda\sigma_3+U,-G^{-1}G_t+ 2i
\lambda(G^{-1}VG-G^{-1}G_{\bar z})\bigg]\nn\\
&=&G^{-1}i\lambda\bigg({1\over2}[S_{z},S_{\bar z}]+2(\al
S)_z\bigg)G. \label{FA3}
\end{eqnarray}
If we set
\begin{eqnarray}{\wt V}=-G^{-1}G_{\bar z}+G^{-1}VG,  \label{wV}
\end{eqnarray}
then (\ref{FA1}) reads
$$
-{\wt V}_z+U_{\bar z}-[{\wt V},\la \sigma_3+U]=0.$$ On the one hand,
from the relation $S=-iG\si_3G^{-1}$ we have
\begin{eqnarray}
 G^{-1}\bigg((i\lambda S)_{\bar z}-V_{z} +[i\lambda S,V]
\bigg)G&=&G^{-1}\bigg(\la(iG\sigma_3G^{-1})_{\bar z}-V_{z} +[\la
G\si_3G^{-1},V] \bigg)G\nn\\
&=&\la [G^{-1}G_{\bar z}-G^{-1}VG,\si_3]-G^{-1}V_zG.\label{UV1}
\end{eqnarray}
On the other hand,  from the relation $G^{-1}G_z=-U$ we also have
\begin{eqnarray}
&& U_{\bar z}- {\widetilde V}_{z}+[\lambda \sigma_3+U,{\widetilde
V}]\nn\\
&=&-(G^{-1}G_z)_{\bar z}-(-G^{-1}G_{\bar z}+G^{-1}VG)_z+[\lambda
\sigma_3-G^{-1}G_z,-G^{-1}G_{\bar z}+G^{-1}VG]\nn\\
&=& -G^{-1}V_zG+\la[G^{-1}G_{\bar z}-G^{-1}VG,\si_3]. \label{UV2}
\end{eqnarray}
Hence, combine (\ref{UV1}) with (\ref{UV2}), we obtain
\begin{eqnarray}
 G^{-1}\bigg((i\lambda S)_{\bar z}-V_{z} +[i\lambda S,V]
\bigg)G=- {\widetilde V}_{z}+
 U_{\bar z}-[{\widetilde V},\lambda \sigma_3+U]
 =\hbox{left-hand-side~of~(\ref{FA1})}, \nn
\end{eqnarray}
which implies that (\ref{FA1}) is automatically satisfied from
(\ref{rV}). We also claim that (\ref{FA2}) is automatically
satisfied. In fact, from the definition (\ref{wV}) of ${\wt V}$, we
have $V=G{\wt V}G^{-1}+G_{\bar z}G^{-1}$ and hence, by the aid of
(\ref{GG}),
$$2V+S_{\bar z}+2\al S=2G{\wt V}G^{-1}.$$
Thereofore
\begin{eqnarray}
&&\hbox{RHS~of~(\ref{FA2})}\nn\\&=& G^{-1}\bigg(-V_t+i\lambda
(V_{\bar z}+S_{\bar z}S+2\al S)_{\bar z}-\bigg[V,i\lambda (V_{\bar
z}+S_{\bar z}S+2\al S)\bigg]
\bigg)G\nn\\
&=&G^{-1}\bigg(-(G{\wt V}G^{-1}+G_{\bar z}G^{-1})_t+2i\la (G{\wt V}G^{-1})_{\bar z}
-[G{\wt V}G^{-1}+G_{\bar z}G^{-1},2i\la G{\wt V}G^{-1}]\bigg)G\nn\\
&=&-{\wt V}_t+(-G^{-1}G_t+2i\la {\wt V})_{\bar z}-[{\wt
V},-G^{-1}G_t+2i\la {\wt V}]\nn\\
&=&\hbox{LHS~of~(\ref{FA2})}.\nn
\end{eqnarray}
So nothing new is obtained from the two identities (\ref{FA1}) and
(\ref{FA2}). Now we come to treat (\ref{FA3}). By using the fact
$\al=-i{\bar p}$, the first equation of (\ref{rq=r}) and the
identities: $S=-iG\si_3G^{-1}$,
$S_z=2iG(U^{\hbox{(off-diag)}}\sigma_3)G^{-1}$ and $S_{\bar
z}=2iG\sigma_3P^{\hbox{(off-diag)}}G^{-1}$ deduced from (\ref{14})
and (\ref{GP}) respectively,  we have
\begin{eqnarray}
&&{{1\over2}[S_{\bar z},S_{z}]+2(\al S)_{z}} \nn\\
&&=G\left(2[U^{\hbox{(off-diag)}}\sigma_3,\sigma_3P^{\hbox{(off-diag)}}]-2
{\bar p}_{z}\sigma_3+4{\bar p} U^{\hbox{(off-diag)}}\sigma_3\right)G^{-1}\nn\\
&&=G\bigg(2p_{\bar z}\sigma_3+4{\bar p}
U^{\hbox{(off-diag)}}\sigma_3\bigg)G^{-1}=G\bigg(2p_{\bar
z}\sigma_3+ [\sigma_3,{H}]\bigg)G^{-1}, \label{dxt}
\end{eqnarray}
where ${H}=\left(\begin{array}{cc}
       0 & -2q{\bar p}\\
       -2{r}{\bar p}& 0
       \end{array}\right)$.
Thus, (\ref{FA3}) is equivalent to holding
\begin{eqnarray}
-U_t+(-G^{-1}G_t+2i\lambda {\wt V})_{z}
-[\lambda\sigma_3+U,-G^{-1}G_t+ 2i\lambda {\wt V}] =2i\lambda
p_{\bar z}\sigma_3 +i\lambda[\sigma_3,{H}], \nn
\end{eqnarray}
or equivalently,
\begin{eqnarray}
-U_t+(-G^{-1}G_t)_{z}+[U,G^{-1}G_t] + \lambda\left(2iU_{\bar
z}-2ip_{\bar z}\sigma_3+
[\sigma_3,G^{-1}G_t]-i[\sigma_3,{H}]\right)=0. \label{EQ}
\end{eqnarray}
Here we have used the identities (\ref{FA1}): $U_{\bar z}-{\wt
V}_{z}+[\lambda \sigma_3+U, {\wt V}]=0$ in the computation. The
vanishing of the coefficients of $\lambda$ and the constant term in
(\ref{EQ}) lead to
\begin{eqnarray}
&&-U_t+(-G^{-1}G_t)_{z}+[U,G^{-1}G_t]=0,\label{term0} \\
&&2U_{\bar z}-2p_{\bar
z}\sigma_3-i[\sigma_3,G^{-1}G_t]-[\sigma_3,{H}]=0. \label{term1}
\end{eqnarray}
Since $G^{-1}G_t\in su(1,1)$, we may set
$G^{-1}G_t=\left(\begin{array}{cc}iu&
\chi\\\bar\chi&-iu\end{array}\right)$ for some real function $u$ and
complex function $\chi$. Then the equation (\ref{term1}) leads to
$\chi=-i(q_{\bar z}+2\bar p q)$ and $\bar \chi=i(r_{\bar z}-2\bar p
r)$ (the compatibility condition $\chi=\bar{\bar \chi}$ is
equivalent to the second equation of (\ref{rq=r})). Thus
\begin{eqnarray}
&&G^{-1}G_t=i\left(\begin{array}{cc}
        u& -q_{\bar z}-2{\bar p}q\\
       r_{\bar z}-2{\bar p}r& -u
       \end{array}\right)=i\left\{\left(u+U_{\bar
z}^{(\hbox{off-diag})}\right)\sigma_3+H\right\}.  \label{019}
\end{eqnarray}
Substituting (\ref{019}) into the equation (\ref{term0}), we finally
get the three  nonlinear Schr\"odinger-type equations given by
(\ref{r8}). This finishes the proof of the theorem. $\Box$

\bigskip

By the way, substituting (\ref{wV}) and (\ref{019}) into (\ref{AG})
and (\ref{017}) respectively, we obtain explicitly
\begin{eqnarray}
 A^G:= {\widetilde V}d{\bar z}+\bigg(\lambda\sigma_3+U\bigg)d{z} +
\bigg\{2i\lambda {\widetilde V}-i(u +U_{\bar
z}^{(\hbox{off-diag})})\sigma_3-i H\bigg\} dt \label{9}
\end{eqnarray}
and
\begin{eqnarray}
{K}^G&:=& \Bigg\{-{\widetilde V}_t+\left(2i\lambda {\widetilde
V}-i(u
+U_{\bar z}^{(\hbox{off-diag})})\sigma_3-i H\right)_{\bar z}\nonumber\\
&&-[{\widetilde V}, 2i\lambda {\widetilde V}-i(u +U_{\bar
z}^{(\hbox{off-diag})})\sigma_3-i H]\Bigg\}d{\bar z}\wedge dt
\nn\\
&& +i\lambda\bigg(2U^{(\hbox{diag})}_{\bar
z}+[\sigma_3,H]\bigg)d{z}\wedge dt, \label{TK}
\end{eqnarray}
where $\lambda$ is the same spectral parameter as in (\ref{r6}),
$U=\left(\begin{array}{cc}
  p & q\\
  {r}& -p
  \end{array}\right)$, $H=\left(\begin{array}{cc}
       0 & -2q{\bar p}\\
       -2{r}{\bar p}& 0
       \end{array}\right)$ and $\widetilde V$ solves the equation
\begin{eqnarray}
U_{\bar z}- {\widetilde V}_{z}+[\lambda \sigma_3+U,{\widetilde
V}]=0. \label{tildeV}
\end{eqnarray}
For the nonlinear Schr\"odinger-type system (\ref{r8},\ref{rq=r}),
as indicated in the proof of Theorem 1, it is a PDE with prescribed
curvature representation:
\begin{eqnarray}
F_{\widetilde A}=d{\widetilde A}-{\widetilde A}\wedge {\widetilde A}
={\widetilde K}, \label{09}
\end{eqnarray}
where $\wt A=A^G$ and $\wt K=K^G$ are given in (\ref{9}) and
(\ref{TK}) respectively.

Next we shall prove that the above gauge transformation from  the
matrix equation (\ref{r5}) of the Schr\"odinger maps from ${\bf
R}^2$ to ${\bf H}^2$  to the nonlinear Schr\"odinger-type system
(\ref{r8},\ref{rq=r}) is in fact reversible.

\begin{Theorem}
For any {\rm $C^2$}-solution $(p,q,r,u)$ to the nonlinear
Schr\"odinger-type system (\ref{r8},\ref{rq=r}), there is a matrix
$C^2$-function $G\in SU(1,1)$ such that $(p,q,r,u)$ is the gauge
transformed to a {\rm $C^3$}-solution $S$ to the equation (\ref{r5})
of the Schr\"odinger maps from ${\bf R}^2$ to ${\bf H}^2$ by $G$.
Moreover, if we require that the gauge matrix $G$ satisfies
$G|_{t=z=\bar z=0}=I$. Then any $C^m$-solution ($m\ge 2$) to the
nonlinear Schr\"odinger-type system (\ref{r8},\ref{rq=r})
corresponds uniquely to a $C^{m+1}$-solution to the equation
(\ref{r5}) of the Schr\"odinger maps from ${\bf R}^2$ to ${\bf H}^2$
and vice versa.
\end{Theorem}

\noindent{\bf Proof}:  Let $(p,q,r,u)$ be a solution to
Eqs.(\ref{r8},\ref{rq=r}). From the proof of Theorem 1, we see that
Eqs.(\ref{r8},\ref{rq=r}) are in fact the integrability condition of
the following linear system:
\begin{eqnarray}
G_{z}=-GU,\quad G_t=Gi\left((u +U^{(\hbox{off-diag})}_{\bar
z})\sigma_3+H\right),\label{10}
\end{eqnarray}
or equivalently,
\begin{eqnarray}
\left\{\begin{array}{c} G_{x_1}=-G\left(\begin{array}{cc}
  i\hbox{Im}{p} & \psi\\
  {\bar \psi}& -i\hbox{Im}{p}
  \end{array}\right)\qquad\\
 G_{x_2}=-G\left(\begin{array}{cc}
  i\hbox{Re}{p} & \phi\\
  {\bar \phi}& -i\hbox{Re}{p}
  \end{array}\right)\qquad\\
 G_t=G\left(\begin{array}{cc}
       iu  & -iq_{\bar z}-2i{\bar p}q\\
       i{r}_{\bar z}-2i{\bar p}{r}& -iu
       \end{array}\right)
       \end{array}\right.
       \label{r10}
\end{eqnarray}
where $\psi=(q+{\bar r})/2$ and $\phi=i(q-{\bar r})/2$. It should be
point out that  Eq.(\ref{rq=r}) implies that the (right) coefficient
matrix in righthand side of the third equation of (\ref{r10}) is
also an $su(1,1)$-matrix. This indicates that the general solutions
$G$ to (\ref{r10}) (i.e. (\ref{10})) belong to the group $SU(1,1)$.
Since the coefficient matrices in (\ref{r10}) are of $C^1$, we let
$G\in SU(1,1)$ be a fundamental $C^2$-smooth solution to (\ref{10})
(or equivalently (\ref{r10})). We shall use it to make the following
gauge transformation for the connection ${\widetilde A}=A^G$ given
in (\ref{9}) with $(p,q,r,u)$ being given above:
\begin{eqnarray}
A&=&(dG) G^{-1} +G {\widetilde A}G^{-1}\nn\\
&=&\bigg(G_{\bar z}G^{-1}+G{\wt V}G^{-1}\bigg)d{\bar
z}+\bigg(G_zG^{-1}+G(\la\sigma_3+U)G^{-1}\bigg)dz\nn\\
&&+ \bigg(G_tG^{-1}+(2i\lambda {\widetilde V}-i(u +U_{\bar
z}^{(\hbox{off-diag})})\sigma_3-i H)G^{-1}\bigg)dt. \label{11}
\end{eqnarray}
We try to show that the 1-form $A$ defined by (\ref{11}) is exactly
the connection of Eq.(\ref{r5}) given by (\ref{r6}) when $S$ and
$\al$ are suitably determined. In fact, substituting the coefficient
$i\lambda S$ of $d{z}$ of (\ref{r6}) into (\ref{11}) and comparing
the coefficients of $\lambda$ of $d{z}$ in the both sides of
(\ref{11}), we obtain
\begin{eqnarray}
G_{z}=-GU,\quad S=-iG\sigma_3G^{-1} \quad ({\rm hence} \quad
S^2=-I). \label{12}
\end{eqnarray}
The first equation of (\ref{12}) is automatically satisfied because
of the first equation of (\ref{10}). The second equation of
(\ref{12}) is regarded as defining $S$. Now, we have to prove that
the coefficients of $d{\bar z}$ and $dt$ of $A$ defined by
(\ref{11}) are respectively the same coefficients of $d{\bar z}$ and
$dt$ of the connection given in (\ref{r6}), that is,
\begin{eqnarray}
 V&=&G_{\bar z}G^{-1}+G{\wt V} G^{-1},  \label{dy} \\
i\lambda\bigg(2V +S_{\bar z} S+2\al S\bigg)
&=&G_tG^{-1}+G\left(2i\lambda {\widetilde
V}-i\bigg((u+U^{(\hbox{off-diag})}_{\bar
z})\sigma_3+H\bigg)\right)G^{-1}. ~~\label{dt}
\end{eqnarray}
Eq.(\ref{dy}) can be regarded as defining $V$ if we can show that
such a $V$  solves Eq.(\ref{rV}), i.e., for the $S$ being given in
(\ref{12}) we must have
\begin{eqnarray}
 (i\lambda S)_{\bar z}-V_{z} +[i\lambda S,V]=0. \label{SV}
\end{eqnarray}
The proof of (\ref{SV}) is a direct computation. Indeed, by using
the expression of $V$ given in (\ref{dy}) and the fact that $G$
fulfills (\ref{r10}) (this equivalent to having (\ref{14}),
(\ref{GP}) and the third equation of (\ref{r10})), we have
\begin{eqnarray}
 (i\lambda S)_{\bar z}-V_{z} +[i\lambda
 S,V]&=&G\bigg(-P_z-[P,U]-{\wt V}_z+[\la \sigma_3+U,{\wt
 V}]\bigg)G^{-1}\nn\\
 &=&
 G\left(U_{\bar z}- {\widetilde V}_{z}+[\lambda \sigma_3+U,{\widetilde
V}] \right)G^{-1}. \nn\label{dy=dy}
\end{eqnarray}
Here have have used the fact: $U_{\bar z}+P_z+[P,U]=0$. Since ${\wt
V}$ satisfies (\ref{tildeV}), this establishes (\ref{SV}).  For
proving (\ref{dt}), since $G$ satisfies the second equation of
(\ref{10}), it is easy to see that (\ref{dt}) is equivalent to
\begin{eqnarray}
2V +S_{\bar z} S+2{\al} S=2G{\wt V}G^{-1}.  \nn
\end{eqnarray}
Since $G{\wt V}G^{-1}=V-G_{\bar z}G^{-1}$ by (\ref{dy}), this
equation can also be written as
\begin{eqnarray}
S_{\bar z} S+2{\al} S=-2G_{\bar z}G^{-1}.  \label{dt=dt}
\end{eqnarray}
Now we take $\al=-i{\bar p}$ which fulfills the requirement of $\al$
in the definition of the connection (\ref{r6}). Substituting
$S=-iG\sigma_3G^{-1}$, $\al=-i{\bar p}$ and applying $G_{\bar
z}=GP$, we have
\begin{eqnarray}
(S_{\bar z} S+2{\al} S)&=&(-iG\sigma_3G^{-1})_{\bar z}
(-iG\sigma_3G^{-1})+2(-i\bar p)(-iG\sigma_3G^{-1})\nn \\
&=&-G{\bar z}G^{-1}+G\sigma_3P\sigma_3G^{-1}-2\bar pG\sigma_3G^{-1}
\nn\\
&=& -G{\bar z}G^{-1}-GPG^{-1}= -2G_{\bar z}G^{-1}. \nn
\end{eqnarray}
This proves (\ref{dt=dt}) and hence (\ref{dt}). Thus we have proved
that the two connections given by (\ref{11}) and (\ref{r6})
respectively are actually the same one when $S=-iG\sigma_3G^{-1}$
and $\al=-i{\bar p}$. What's the remainder for us to do is to prove
that the curvature formula
\begin{eqnarray}
F_A=K =G{\widetilde K}G^{-1}=GF_{\widetilde A}G^{-1} \label{curv}
\end{eqnarray}
under the gauge transformation is satisfied too, where $K$ is given
by (\ref{rK}) and ${\widetilde K}=K^G$ is given by (\ref{TK}). In
fact, on the one hand, we see that
\begin{eqnarray}
G{\wt K}G^{-1}&=& G\Bigg(\bigg\{-{\widetilde V}_t+\left(2i\lambda
{\widetilde V}-i(u+U^{(\hbox{off-diag})}_{\bar z})
\sigma_3-iH\right)_{\bar z}\nonumber\\
&&-[{\widetilde V}, 2i\lambda {\widetilde
V}-i(u+U^{(\hbox{off-diag})}_{\bar
z})\sigma_3-i H]\bigg\}d{\bar z}\wedge dt\nn\\
&&+\lambda(2U^{(\hbox{diag})}_{\bar z}+[\sigma_3,H])d{z}\wedge
dt\Bigg)G^{-1}. \label{GKG}
\end{eqnarray}
On the other hand, by using (\ref{dy}, \ref{dt=dt}) and (\ref{10}),
it is a straightforward calculation that the coefficient of $d{\bar
z}\wedge dt$ in $K$ given by (\ref{rK}) is
\begin{eqnarray}
&&-V_t+i\lambda (2V +S_{\bar z} S+2{\al} S)_{\bar z}-[V,i\lambda(2V
+S_{\bar z} S+2{\al} S)]\nonumber \\
&=&-V_t+2i\lambda(G{\wt V}G^{-1})_{\bar z}-[V,2i\lambda G{\wt
V}G^{-1}]\nn \\
&=&-G_{t{\bar z}}G^{-1}+G_{\bar z}GG_tG^{-1}-G_t{\wt V}G^{-1}-G{\wt
V}_tG^{-1}+G{\wt V}G^{-1}G_tG^{-1}+2i\lambda G_{\bar z}{\wt V}G^{-1}\nn\\
&&+2i\lambda G{\wt V}_{\bar z}G^{-1}-2i\lambda G{\wt V}G^{-1}G_{\bar
z}G^{-1}-2i\lambda(G_{\bar z}{\wt V}G^{-1}-G{\wt
V}G^{-1}G_{\bar z}G^{-1})\nn\\
&=&G\Bigg\{-{\widetilde V}_t+\left(2i\lambda {\widetilde
V}-i(u+U^{(\hbox{off-diag})}_{\bar z})\sigma_3
-iH\right)_{\bar z}\nonumber\\
&&-[{\widetilde V}, 2i\lambda {\widetilde
V}-i(u+U^{(\hbox{off-diag})}_{\bar z})\sigma_3-i H]\Bigg\}G^{-1}.
\label{dydt}
\end{eqnarray}
(\ref{GKG}) and (\ref{dydt}) indicate that the two coefficients of
$d{\bar z}\wedge dt$ in the both sides of (\ref{curv}) are the same
one. Meanwhile, the coefficient of $dz\wedge dt$ of $K$ is
\begin{eqnarray}
i\lambda\bigg({1\over2}[S_{z},S_{\bar z}]+2(\al
S)_{z}\bigg)=i\lambda G \bigg(2U^{(\hbox{diag})}_{\bar
z}+[\sigma_3,H]\bigg)G^{-1}, \label{dxdt}
\end{eqnarray}
Here we have used the identity: ${1\over2}[S_{\bar z},S_{z}]+2(\al
S)_{z}=G(2p_{\bar z}\sigma_3+ [\sigma_3,{H}])G^{-1}$, where
${H}=\left(\begin{array}{cc}
       0 & -2q{\bar p}\\
       -2{r}{\bar p}& 0
       \end{array}\right)$
as before, which is deduced from the same argument in getting
(\ref{dxt}). (\ref{GKG}) and (\ref{dxdt}) imply that the two
coefficients of $d{z}\wedge dt$ in the both sides of (\ref{curv})
are also the same one. Thus we have proved the desired identity
(\ref{curv}), which implies the holding of the prescribed curvature
representation (\ref{r7}) by the gauge theory. This indicates that
$S$ defined by the second equation of (\ref{12}) from the solution
$(p,q,r,u)$ to (\ref{r8},\ref{rq=r}) satisfies the matrix (\ref{r5})
and hence the equation of the Schr\"odinger maps from ${\bf R}^2$ to
${\bf H}^2$.

Since $G$ is a solution to the linear first-order differential
system (\ref{r10}), it is well-known from linear theory of
differential equations that such a $G$ is unique if we propose the
initial condition $G|_{t=z=\bar z=0}=I$ on $G$. Under this
circumstance, we see that a solution $(p,q,r,u)$ to
Eqs.(\ref{r8},\ref{rq=r}) corresponds uniquely to a solution $S$ to
Eq.(\ref{r5}) by the gauge transformation and vice versa.
Furthermore, because of the relation
$S_z=2iGU^{(\hbox{off-diag})}\sigma_3G^{-1}$ and $S_{\bar
z}=-2iGP^{(\hbox{off-diag})}\sigma_3G^{-1}$ deduced from (\ref{12}),
the remainder part of the theorem is obviously true. $\Box$

\begin{Remark}
Theorem 1,2 reveal a closed relation between $C^m$-solutions $S$ to
the Schr\"odinger maps from ${\bf R}^2$ to ${\bf H}^2$ and
$C^{m-1}$-solutions $(p,q,r,u)$ to the nonlinear Schr\"odinger-type
system (\ref{r8},\ref{rq=r}). As we shall see in the next section,
the nonlinear Schr\"odinger-type system (\ref{r8},\ref{rq=r}) will
play an essential role for us to show the existence of blowing-up
$C^{\infty}$-solutions to the Schr\"odinger maps from ${\bf R}^2$ to
${\bf H}^2$.
\end{Remark}

Let us introduce the polar coordinates $(\rho,\theta)$ of ${\bf
R}^2$, that is, $x_1=\rho\cos\theta, x_2=\rho\sin\theta$. Thus we
have
$${{\partial}\over{\partial z}}=
{{e^{-i\theta}}}\left({{\partial}\over{\partial \rho}}-{i\over
\rho}{{\partial}\over{\partial \theta}}\right),\quad
{{\partial}\over{\partial {\bar z}}}=
{{e^{i\theta}}}\left({{\partial}\over{\partial \rho}}+{i\over
\rho}{{\partial}\over{\partial \theta}}\right).$$ We would like to
find following ansatz solutions to (\ref{r8},\ref{rq=r}):
\begin{eqnarray}
p=0, q=e^{-i\theta}Q(t,\rho), r=e^{-i\theta}{\bar Q}(t,\rho),
u=|Q|^2(t,\rho)-2\int_{\rho}^{\infty}{{|Q|^2(t,\tau)}\over{\tau}}d\tau\label{ansatz}
\end{eqnarray}
for some suitable functions $Q(t,\rho)$. One may verify that
(\ref{rq=r}) and the third equation of (\ref{r8}) are satisfied
automatically (since $
\partial^{-1}_z\partial_{\bar z}(rq)=\partial^{-1}_z\partial_{\bar
z}(e^{-2i\theta}|Q|^2(t,\rho))=|Q|^2(t,\rho)-2\int_{\rho}^{\infty}
{{|Q|^2(t,\tau)}\over{\tau}}d\tau $)  and the first and second
equation of (\ref{r8}) lead $Q$ to satisfy the nonlinear
Schr\"odinger differential-integral equation (\ref{QQ}) presented in
Introduction. This equation was deduced in \cite{CSU} by the
(generalized) Hasimoto transformation. Furthermore, under this
circumstance the linear system (\ref{r10}) is reduced to
\begin{eqnarray}
\left\{\begin{array}{c}G_{\rho}=-G\left(\begin{array}{cc}0& Q(t,\rho)\\
\bar
Q(t,\rho)&0\end{array}\right),~~~~~~~~~~~~~~~~~~~~~~~~~~~~~~~~~~~
~~~~~~~~~~~~~~~~\\
G_t=Gi\left(\begin{array}{cc} |Q|^2(t,\rho)-2\int_{\rho}^{\infty}
{{|Q|^2(t,\tau)}\over{\tau}}d\tau & -Q_{\rho}(t,\rho)-Q(t,\rho)/\rho\\
\bar Q_{\rho}(t,\rho)+\bar Q(t,\rho)/\rho
&-|Q|^2(t,\rho)+2\int_{\rho}^{\infty}
{{|Q|^2(t,\tau)}\over{\tau}}d\tau\end{array}\right).\end{array}\right.\label{Grho}
\end{eqnarray}
So the gauge matrix function $G(z,\bar z,t)=G(t,\rho)$ is radial.
The corresponding matrix $S$ given by the second equation (\ref{12})
is exactly $-iG(t,\rho)\sigma_3G^{-1}(t,\rho)$ and hence is a
radially symmetric solution to Eq.(\ref{r5}) of the Schr\"odinger
flow from ${\bf R}^2$ to the hyperbolic 2-space ${\bf H}^2$.

\section * {\S 3. Applications}

We will follow the basic conventional notations for Sobolev spaces
 $W^{m,\sigma}({\bf R}^n)$ ($\sigma\ge2$) of real
or complex-valued functions or spaces $C^k({\bf R}^n)$ of continuous
differential functions up to order $k$ on ${\bf R}^n$ for $n\ge2$
and norms $||\cdot||_{W^{m,\sigma}({\bf R}^n)}$ or
$||\cdot||_{C^k({\bf R}^n)}$ used in \cite{GiT}

The gauged nonlinear Schr\"odinger-type system (\ref{r8},\ref{rq=r})
looks much complicated and it seems seldom useful. However, the
nonlinear Schr\"odinger-type system (\ref{r8},\ref{rq=r}) is in fact
another mathematical description of the Schr\"odinger maps from
${\bf R}^2$ to the hyperbolic 2-space ${\bf H}^2$ and should be
useful. In this section, as applications, we shall use the gauged
nonlinear Schr\"odinger-type system (\ref{r8},\ref{rq=r}) to explore
some geometric property of harmonic maps from ${\bf R}^2$ to ${\bf
H}^2$ and especially the blow-up property of the Schr\"odinger maps
from ${\bf R}^2$ to the hyperbolic 2-space ${\bf H}^2$.

\bigskip
The first application is to study whether there exists a nontrivial
smooth radially symmetric harmonic map from the whole plane ${\bf
R}^2$ to ${\bf H}^2$. If there does exist a such harmonic map, then
we may have a nontrivial radial solitonary wave solution $(p,q,r,u)$
to (\ref{r8},\ref{rq=r}) from it. Unfortunately, there is no such a
nontrivial radial symmetric smooth harmonic maps from ${\bf R}^2$ to
${\bf H}^2$ defined on the whole plane ${\bf R}^2$.

\begin{Theorem}
There does not exist a nontrivial radially symmetric smooth harmonic
map from the whole plane ${\bf R}^2$ to ${\bf H}^2$.
\end{Theorem}
\nd{\bf Proof}. If there exists a nontrivial radially symmetric
smooth harmonic map from ${\bf R}^2$ to ${\bf H}^2$, then we denote
it by $S=S(\rho)$, where $(\rho,\theta)$ is the polar coordinate of
the plane. From the appendix at the end of the paper, we see that
\begin{eqnarray}
G={1\over{\sqrt{2(1-s_3)}}}(S-i\sigma_3){\rm
diag}(\gamma,{\bar\gamma}),\nn \label{gammas}
\end{eqnarray}
is a smooth $SU(1,1)$-matrix function satisfying
$\sigma_3=iG^{-1}SG$, where $\gamma$  is a complex function of
$\rho$ and $t$ with $|\gamma|=1$. Since
${1\over{\sqrt{2(1-s_3)}}}(S-i\sigma_3)$ depends only on the radial
variable $\rho$ and is nontrivial, we may choose $\gamma$ such that
\begin{eqnarray}
G^{-1}G_{\rho}=- \left(\begin{array}{cc}
       0 & q(t,\rho)\\
       \bar q(t,\rho)& 0
       \end{array}\right),\quad
G^{-1}G_t=\left(\begin{array}{cc}
  \gamma^{-1}\gamma_t & 0\\
  0& {\bar\gamma}^{-1}{\bar\gamma}_t
  \end{array}\right), \label{us}
\end{eqnarray}
for some nonzero function $q(t,\rho)$ depending on $\rho$ and
possibly on $t$. On the other hand, from the proof of Theorem 1 we
know such $G$ should satisfy (\ref{Grho}) for some smooth
(nontrivial) $Q(t,\rho)$ too. Notice that the two matrices $G$
satisfying respectively (\ref{us}) and (\ref{Grho}) are at most up
to a (smooth) diagonal matrix ${\rm diag}(\gamma_0,{\bar\gamma_0})$
with $|\gamma_0|=1$. Thus, compare (\ref{Grho}) with (\ref{us}), we
have
\begin{eqnarray}
Q(t,\rho)=\gamma_0^{-1}q(t,\rho)\bar{\gamma_0},\quad
Q_{\rho}+Q/\rho=0.
\end{eqnarray}
Hence $Q(t,\rho)=\frac{C(t)}{\rho}$ for some nonzero constant $C(t)$
depending only on $t$. Since $Q(t,\rho)$ has singular at the origin
except $C(t)=0$, which is a contradiction to that $Q(t,\rho)$ and
hence $q(t,\rho)$ is a nontrivial smooth function of $t$ and $\rho$.
This contradiction implies that $S=S(\rho)$ is not nontrivial smooth
on the whole plane, which proves the theorem.~~$\Box$

\bigskip
The next application is to use the gauged nonlinear
Schr\"odinger-type system (\ref{r8},\ref{rq=r}) to construct a class
of blowing-up solutions of the Schr\"odinger maps to the hyperbolic
2-space ${\bf H}^2$ (\ref{r5}). This is an interesting problem in
the study of Schr\"odinger maps which was proposed by Y.W. Ding in
\cite{Di} as a unsolved problem. In order to study such a problem,
let's first give a brief description of some blowing-up results for
nonlinear Schr\"odinger equations. It is well-known that, a
nonlinear Schr\"odinger equation with the critical Sobolev
exponential $\si=1+4/n$ in ${\bf R}^n$: $iq_t+\Delta
q+|q|^{\si-1}q=0$ admits a so-called conformal invariance (see, for
example, \cite{Wei,bib:wein,SuS}). That says, the transformation
from a solution $q$  to
$$
{\wt
q}(t,x)={{e^{i{{b|x|^2}\over{4(a+bt)}}}}\over{(a+bt)^{n/2}}}q(T,X),
$$
is invariant, where $X={{x}\over{a+bt}}$ with
$x=(x_1,\cdots,x_n)\in{\bf R}^n$, $T={{c+dt}\over{a+bt}}$ and
$\left(\begin{array}{cc}
              a & b\\
              c& d\end{array} \right)\in SL_2({\bf R})$ (i.e.
              $a, b, c, d$ are real and $ad-bc=1$), i.e.,
${\wt q}(t,x)$ is still a solution to the same equation. Weinstein
constructed in \cite{Wei} blowing up solutions to the nonlinear
Schr\"odinger equation with the critical exponential from solitonary
wave solutions (see \cite{Li,Stra} for the existence of such
solitonary wave solutions) by using this conformal invariant
property. Now for our present nonlinear Schr\"odinger-type system
(\ref{r8},\ref{rq=r}), the following lemma is very crucial. This
fact is somewhat difficult to take note of, but it arises
heuristically from the conformal invariant property of nonlinear
Schr\"odinger equations mentioned above.

\begin{Lemma}
Let $({p},{q},{r},{u})$ be given as follows:
\begin{eqnarray}
p={p}(t,z,\bar z)&=&
0\nn\\
q={q}(t,z,\bar z)&=&{{e^{i{{bz\bar z}\over{2(a+bt)}}}}\over{a+bt}}
\al \bar z\nn\\
r={r}(t,z,\bar z)&=&
{{e^{-i{{bz\bar z}\over{2(a+bt)}}}}\over{(a+bt)}}\al \bar z\nn\\
u={u}(t,z,\bar z)&=&{{2\al^2z{\bar z} }\over{(a+bt)^2}},\nn
\label{45qq}
\end{eqnarray}
where $a$, $b$ and $\al$ are real numbers. Then $({p},{q},{r},{u})$
is a solution to the gauged nonlinear Schr\"odinger system
(\ref{r8},\ref{rq=r}) if and only if  $\al^2=\frac{b^2}{16}$.
\end{Lemma}

\nd {\bf Proof}. The proof is just a direct verification one by one:
\begin{eqnarray}
-\bar r_z+q_{\bar z}+2(p\bar r+\bar p q)=-\left({{e^{i{{bz\bar
z}\over{2(a+bt)}}}}\over{a+bt}} \al z\right)_z+\left({{e^{i{{bz\bar
z}\over{2(a+bt)}}}}\over{a+bt}} \al \bar z\right)_{\bar z} =0,\nn
\end{eqnarray}
\begin{eqnarray}
\bar p_z+p_{\bar z}+|q|^2-|r|^2=\left|\frac{{e^{i{{bz\bar
z}\over{2(a+bt)}}}}\al\bar
z}{a+bt}\right|^2-\left|\frac{{e^{-i{{bz\bar
z}\over{2(a+bt)}}}}\al\bar z}{a+bt}\right|^2=0,\nn
\end{eqnarray}
and
\begin{eqnarray}
&&i{q}_t+{q}_{z\bar z}-2{u}{q}+2({\overline{p}}{q})_z-2{p} q_{\bar
z}-4|{p}|^2{q}\nn \\
 &&\quad=i\left({{e^{i{{bz\bar
z}\over{2(a+bt)}}}}\over{a+bt}} \al \bar
z\right)_t+\left({{e^{i{{bz\bar z}\over{2(a+bt)}}}}\over{a+bt}} \al
\bar z\right)_{z\bar z}-{{4\al^2z{\bar z} }\over{(a+bt)^2}}
\left({{e^{i{{bz\bar
z}\over{2(a+bt)}}}}\over{a+bt}} \al \bar z\right)\nn\\
&&\quad=\left(\frac{b^2}{2}-\frac{b^2}{4}-4\al^2\right){{e^{i{{bz\bar
z}\over{2(a+bt)}}}}\over{(a+bt)^3}} \al \bar z; \nn\\
&&i{r}_t-{r}_{z\bar z}+2{u}{r}+2({\overline{p}}{r})_z-2{p}{r}_{\bar
z}+4|{p}|^2{r}\nn \\
 &&\quad= i\left({{e^{-i{{bz\bar
z}\over{2(a+bt)}}}}\over{a+bt}} \al \bar
z\right)_t-\left({{e^{-i{{bz\bar z}\over{2(a+bt)}}}}\over{a+bt}} \al
\bar z\right)_{z\bar z}+{{4\al^2z{\bar z} }\over{(a+bt)^2}}
\left({{e^{-i{{bz\bar
z}\over{2(a+bt)}}}}\over{a+bt}} \al \bar z\right)\nn\\
&&\quad=\left(-\frac{b^2}{2}+\frac{b^2}{4}+4\al^2\right){{e^{-i{{bz\bar
z}\over{2(a+bt)}}}}\over{(a+bt)^3}} \al \bar z;\nn\\
&&ip_t+(qr)_{\bar z}-u_z=\left(\frac{\al^2(\bar
z)^2}{(a+bt)^2}\right)_{\bar z}-\left({{2\al^2z{\bar z}
}\over{(a+bt)^2}}\right)_z=0.\nn
\end{eqnarray}
Since $\al^2=\frac{b^2}{16}$, the lemma is thus proved. $\Box$

\bigskip
It is obvious that solutions $(p,q,r,u)$ constructed in Lemma 1 blow
up in finite time (as $t\to -\frac{a}{b}$) at every point in the
plane except the origin. We now prove the existence of blowing-up
solutions to Eq.(\ref{r5}) from these solutions.

\begin{Theorem}
For real numbers $a$ and $b$ with $a>0$ and $b<0$, there exists a
$C^{\infty}$-solution $S(t,z,\bar z)$ to the equation (\ref{r5}) of
the Schr\"odinger maps from ${\bf R}^2$ to the hyperbolic 2-space
${\bf H}^2$ such that the absolute of the gradient $|\nabla S|$
blows up in finite time at every point in ${\bf R}^2$ except the
origin.
\end{Theorem}
\nd{\bf Proof}. For real numbers $a$ and $b$ with $a>0$ and
              $b<0$,
we may construct the explicit solution $({p},{q},{r},{u})$ to the
gauged nonlinear Schr\"odinger system (\ref{r8},\ref{rq=r}) by
$$
p=0,~ q={{e^{i{{bz\bar z}\over{2(a+bt)}}}}\over{a+bt}} \al \bar z,~
r={{e^{-i{{bz\bar z}\over{2(a+bt)}}}}\over{(a+bt)}}\al \bar
z,~u={{2\al^2z{\bar z} }\over{(a+bt)^2}}$$ from Lemma 1, where
$\al^2=\frac{b^2}{16}$. By Theorem 2  there is a smooth solution
$S(z,\bar z,t)$ to the Eq.(\ref{r5}) which is the $SU(1,1)$-gauge
equivalent to the above solution $({p},{q},{r},{u})$ of
(\ref{r8},\ref{rq=r}) by a gauge matrix $G$. On the one hand, from
the formulas: $G_{z}=-GU$, $G_{\bar z}=GP$ displayed in Theorem 1 or
2, where $U=\left(\begin{array}{cc}
  {p} & {{q}}\\
  { r}& -{ p}
  \end{array}\right)$, $P=\left(\begin{array}{cc}
  \overline{{ p}} & \overline{{{ r}}}\\
  \overline{{ q}}& -\overline{{ p}}
  \end{array}\right)$ and $G\in SU(1,1)$, we have $$S_{z}=2iG
U^{(\hbox{off-diag})}\sigma_3G^{-1},~~~S_{\bar z}=-2iG
P^{(\hbox{off-diag})}\sigma_3G^{-1}.$$ On the other hand, from the
appendix at the end of this paper we have that such $G$ can be
expressed as
$$
G={1\over{\sqrt{2(1-s_3)}}}(S-i\sigma_3){\rm
diag}(\gamma,{\bar\gamma}), $$ where $\gamma$ is a complex function
of $z$, $\bar z$ and $t$ with $|\gamma|=1$. Hence, the above two
matrix identities read
\begin{eqnarray}
(S-i\sigma_3){\rm diag}(\gamma,{\bar\gamma})\left(\begin{array}{cc}
  0 & -q\\
  r& 0
  \end{array}\right)=\frac{1}{2i}S_{z}(S-i\sigma_3){\rm
  diag}(\gamma,{\bar\gamma})
   \nn\end{eqnarray}
   and
\begin{eqnarray}
(S-i\sigma_3){\rm diag}(\gamma,{\bar\gamma})\left(\begin{array}{cc}
  0 & \bar r\\
  -\bar q& 0
  \end{array}\right)=\frac{1}{2i}S_{\bar z}(S-i\sigma_3){\rm
  diag}(\gamma,{\bar\gamma}).
   \nn\end{eqnarray}
Equaling the corresponding (row 1 and column 2) entries in the both
sides of the matrix identities, we obtain
\begin{eqnarray}
-i(s_3-1)q\ga&=&\frac{1}{2}\bigg[(s_3)_z(s_1+is_2)-
((s_1)_z+i(s_2)_z)(s_3-1)\bigg]\bar\ga\nn\\
i(s_3-1)\bar r\ga&=&\frac{1}{2}\bigg[(s_3)_{\bar
z}(s_1+is_2)-((s_1)_{\bar z}+i(s_2)_{\bar
z})(s_3-1)\bigg]\bar\ga.\nn
  \end{eqnarray}
Hence
\begin{eqnarray}
|q|+|r|\le
  |(s_1)_{z}|+|(s_1)_{\bar z}|+|(s_2)_{z}|+
  |(s_2)_{\bar z}|+|(s_3)_{z}|+|(s_2)_{\bar z}|:=|\nabla S|.\nn
  \end{eqnarray}
Here we have used the facts: $|\ga|=1$, $|s_1/(1-s_3)|\le 1$ and
$|s_2/(1-s_3)|\le 1$  because of $s_3^2=1+s_1^2+s_2^2$ and $s_3<0$.
Since $q={{e^{i{{bz\bar z}\over{2(a+bt)}}}}\over{a+bt}} \al \bar z$
and $r={{e^{-i{{bz\bar z}\over{2(a+bt)}}}}\over{a+bt}} \al \bar z$
given above, we see that $|q|+|r|=2\frac{|\al|\cdot|z|}{(a+bt)}\to
+\infty$ as $t\to -\frac{a}{b}$ except $z=0$. This implies that
$|\nabla S|$ blows up at every point in the plane as $t\to
-\frac{a}{b}$ except $z=0$. The proof of the theorem is completed.
$\Box$

\begin{Remark}
The collection of solutions $S(t,z,\bar z)$ constructed in Theorem 4
composes a class of blowing-up $C^{\infty}$-solutions to the
Schr\"odinger maps from ${\bf R}^2$ to the hyperbolic 2-space ${\bf
H}^2$. Since these solutions are constructed from explicit solutions
$(p,q,r,u)$ to the gauged nonlinear Schr\"odinger-type system
(\ref{r8},\ref{rq=r}) in Lemma 1, they are also called explicit
blowing-up solutions to the Schr\"odinger maps from ${\bf R}^2$ to
the hyperbolic 2-space ${\bf H}^2$ in this paper.
\end{Remark}

\section *{\S 4. Conclusion}
In this paper, we have proved the existence of a class of smooth
blowing-up solutions to the Schr\"odinger maps from ${\bf R}^2$ to
the hyperbolic 2-space ${\bf H}^2$ with aid of its gauged equivalent
nonlinear Schr\"odinger-type system (\ref{r8},\ref{rq=r}). This sets
a concrete example of the existence of blowing-up smooth solutions
to Schr\"odinger maps, which answers a problem proposed in \cite{Di}
for Schr\"odinger flows. It also indicates that the introduction of
the gauged nonlinear Schr\"odinger-type system (\ref{r8},\ref{rq=r})
is much more effective than that of the known Eq.(\ref{QQ}) and we
believe that this system will be much useful in exploring some
deeper dynamical properties of the Schr\"odinger maps from ${\bf
R}^2$ to ${\bf H}^2$. By the way, in a completely similar way, one
may transform the $1+2$ dimensional Landau-Lifshitz equation: ${\bf
s}_t={\bf s}{\times}\left({{\partial^2}\over{\partial x_1^2}}+
{{\partial^2}\over{\partial x_2^2}}\right){\bf s}$, where
$s=(s_1,s_2,s_3)$ $\in S^2\hookrightarrow {\bf R}^3$, into the
following $1+2$ dimensional nonlinear Schr\"odinger-type system
\begin{eqnarray}
\left\{\begin{array}{c} iq_t-q_{z{\bar z}}+2u
q-2({\bar p}q)_z+2pq_{\bar z}+4|p|^2q=0\\
ir_t+r_{z{\bar z}}-2u r-2({\bar p}r)_z+2pr_{\bar z}-4|p|^2r=0\\
ip_t=(qr)_{\bar z}-u_z\qquad\qquad\qquad\qquad\qquad\qquad \\{\bar
p}_z+p_{\bar z}= |q|^2-|r|^2~\qquad\qquad \qquad\qquad\qquad  \\
{\bar r}_z+q_{\bar z}=2(p{\bar r}-{\bar p}q),\qquad\qquad
\qquad\qquad\qquad
\end{array} \right.    \nn
\end{eqnarray}
by an $SU(2)$-gauge matrix.  From this system we may also prove that
there are no nontrivial radially symmetric smooth harmonic maps from
the whole ${\bf R}^2$ to the 2-sphere ${\bf S}^2$. But we are unable
to establish a similar Lemma 1 for the present system and hence not
to explicitly construct blowing-up smooth solutions to the
Landau-Lifshitz equation.

Finally, we remark that the blowing-up smooth solutions to the
Schr\"odinger maps from ${\bf R}^2$ to ${\bf H}^2$ constructed in
this paper do not belong to any Sobolev space ${W^{m,\sigma}({\bf
R}^2)}$ (for a fixed $t$). It is very interesting that whether there
exist blowing-up ${W^{m,\sigma}({\bf R}^2)}$-solutions to
Schr\"odinger maps, especially to the Landau-Lifshitz equations.
However, this problem is still unknown.

\section * {\bf Acknowledgement}

This work is partially supported by NNFSC (10531090).

\section*{Appendix}
This appendix gives a proof of the following fact used in the proofs
of Theorem 1 and 4.

\begin{Proposition}
For a given matrix $S=S(t,z,\bar z)=\left(\begin{array}{cc}
       is_3 &s _1+is_2\\
       s_1-is_2 & -is_3
       \end{array}\right)\in su(1,1)$,
where $(s_1)^2+(s_2)^2-(s_3)^2=-1$ with $s_3<0$, the general
$SU(1,1)$-solutions $G=G(t,z,\bar z)$ to the matrix equation
$\sigma_3=iG^{-1}SG$ are of the from:
\begin{eqnarray}
G={1\over{\sqrt{2(1-s_3)}}}(S-i\sigma_3){\rm
diag}(\gamma,{\bar\gamma}), \nn
\end{eqnarray}
where $\gamma$  is a complex function of $z$, $\bar z$ and $t$ with
$|\gamma|=1$.
\end{Proposition}
\nd{\bf Proof}. Since a matrix in $SU(1,1)$ is of the form
$\left(\begin{array}{cc}
       A & B\\
       \bar B& {\bar A}
       \end{array}\right)$, where $A$ and $B$ are some complex
       numbers with $|A|^2-|B|^2=1$, we may assume that $G=\left(\begin{array}{cc}
       \nu & \chi\\
       \bar \chi& {\bar \nu}
       \end{array}\right)$, where $\nu$ and $\chi$ are unknown complex
functions of $z$, $\bar z$ and $t$ with $|\nu|^2-|\chi|^2=1$.
Substituting this matrix expression into the matrix equation
$\sigma_3=iG^{-1}SG$, we see that the matrix equation is equivalent
to
\begin{eqnarray}
\left\{\begin{array}{c} -(s_3+1)\nu+i(s_1+is_2)\bar\chi=0\\
i(s_1-is_2)\nu+(s_3-1)\bar\chi=0.\end{array}\right.\label{lsystem}
\end{eqnarray}
Since $(s_1)^2+(s_2)^2-(s_3)^2=-1$ and hence the rank of the
coefficient matrix of the linear system (\ref{lsystem}) is 1, we
have from linear algebra that the linear system (\ref{lsystem}) has
nontrivial solutions and all the solutions compose an 1-dimensional
linear (complex) space. Now it is a direct verification that a pair
of $\nu=\nu_0=i(s_3-1)/\sqrt{2(1-s_3)}$ and $\bar \chi=\bar
\chi_0=(s_1-is_2)/\sqrt{2(1-s_3)}$ solves the linear system
(\ref{lsystem}) and satisfies $|\nu_0|^2-|\chi_0|^2=1$. So
$(\nu_0,\chi_0)$ can be regarded as a base of the linear (complex)
space of solutions to (\ref{lsystem}). Hence the general solutions
to (\ref{lsystem}) are
$$(\nu,\bar\chi)=\gamma(\nu_0,\bar\chi_0)$$
for a complex function $\ga=\gamma(t,z,\bar z)$ and the restriction
of $|\nu|^2-|\chi|^2=1$ implies $|\ga|=1$. Return to matrix $G$, we
thus obtain that the general $SU(1,1)$-solutions to
$\sigma_3=iG^{-1}SG$ are
\begin{eqnarray}
G={1\over{\sqrt{2(1-s_3)}}}(S-i\sigma_3){\rm
diag}(\gamma,{\bar\gamma}), \nn
\end{eqnarray}
where $\gamma$  is a complex function  with $|\gamma|=1$. This
completes the proof of the proposition.~~$\Box$

\end{document}